\documentclass[a4paper,11pt]{article}
\usepackage{authblk}

\usepackage{amssymb}
\usepackage{graphicx}

\usepackage{geometry}
\usepackage{url}

\usepackage{color}
\usepackage{amssymb} 
\usepackage{amsmath}
\usepackage{tikz}
\usetikzlibrary{calc,arrows}

\renewcommand{\vec}[1]{\mathbf{#1}}

\newcommand{\mat}[1]{\mathbf{#1}}
\DeclareMathAlphabet\mathbfcal{OMS}{cmsy}{b}{n} 
\newcommand{\ten}[1]{\mathbfcal{#1}}
\usepackage{relsize} 
\renewcommand{\top}{\mathsmaller{T}}

\newcommand{\llbrack}{\left[\!\left[}
\newcommand{\rrbrack}{\right]\!\right]}

\newcommand{\eg}{\emph{e.g., }} 
\newcommand{\ie}{\emph{i.e., }} 

\newcommand{\citep}[1]{\cite{#1}}
\newcommand{\citet}[1]{\cite{#1}}

\begin{document}


\title{Decoupling multivariate functions\\using second-order information and tensors}



\author{Philippe Dreesen}
\author{Jeroen {De Geeter}}
\author{Mariya Ishteva}

\affil{Vrije Universiteit Brussel (VUB), Dept.~VUB-ELEC, Brussels, Belgium}



\date{}

\maketitle

\begin{abstract}
The power of multivariate functions is their ability to model a wide variety of phenomena, but have the disadvantages that they lack an intuitive or interpretable representation, and often require a (very) large number of parameters. 
We study decoupled representations of multivariate vector functions, which are linear combinations of univariate functions in linear combinations of the input variables. 
This model structure provides a description with fewer parameters, and reveals the internal workings in a simpler way, as the nonlinearities are one-to-one functions. 
In earlier work, a tensor-based method was developed for performing this decomposition by using first-order derivative information. 
In this article, we generalize this method and study how the use of second-order derivative information can be incorporated. 
By doing this, we are able to push the method towards more involved configurations, while preserving uniqueness of the underlying tensor decompositions.
Furthermore, even for some non-identifiable structures, the method seems to return a valid decoupled representation. 
These results are a step towards more general data-driven and noise-robust tensor-based  framework for computing decoupled function representations.
\end{abstract}

\section{Introduction}
\subsection{Towards interpretability of nonlinear models}
Nonlinear models are used in a wide variety of science and engineering fields, such as data analytics, signal processing, system identification, and control engineering. 
While nonlinear models are able to capture wild nonlinear effects, this often comes at the cost of high parametric complexity, and a lack of `model interpretability'.

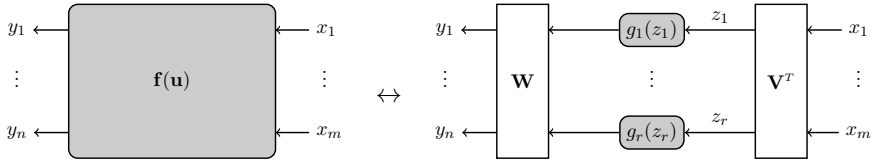
\begin{figure*}[!htb]
\begin{center}
\tikz \node[scale=0.68]{\begin{tikzpicture}
\node (y1) at (0,3) {$y_1$};
\node at (0,2.15) {$\vdots$};
\node (yn) at (0,1) {$y_n$};
\draw [thick,fill=black!20, rounded corners=5pt] (1,0.5) rectangle (5,3.5); \node (F) at (3,2) {$\vec{f}(\vec{u})$};
\draw [<-, thick, label=] (5,3) -- (5.65,3) node[right] {$x_1$};
\node at (5.95,2.15) {$\vdots$};
\draw [<-, thick] (5,1) -- (5.65,1) node[right] {$x_m$};
\draw [<-, thick] (y1) -- (1,3);
\draw [<-, thick] (yn) -- (1,1);
\end{tikzpicture}
\quad
\raisebox{6\height}{\Large$\leftrightarrow$}
\quad
\begin{tikzpicture}
\node (y1) at (0,3) {$y_1$};
\node at (0,2.15) {$\vdots$};
\node (yn) at (0,1) {$y_n$};
\draw [thick] (1,0.5) rectangle (2,3.5); \node (L) at (1.5,2) {$\mat{W}$};
\draw [<-, thick] (y1) -- (1,3);
\draw [<-, thick] (yn) -- (1,1);
\node [shape=rectangle,draw,thick,fill=black!20,rounded corners=5pt] (g1) at (4,3) {$g_1(z_1)$};
\draw [<-, thick, label=] (2,3) -- (g1) node[above,midway] { };
\node at (4,2.15) {$\vdots$};
\node [shape=rectangle,draw,thick,fill=black!20,rounded corners=5pt] (gr) at (4,1) {$g_r(z_r)$};
\draw [<-, thick] (2,1) -- (gr) node[above,midway] { };
\draw [thick] (6,0.5) rectangle (7,3.5); \node (R) at (6.5,2) {$\mat{V}^\top$};
\draw [<-, thick] (g1) -- (6,3) node[above,midway] {$z_1$};
\draw [<-, thick] (gr) -- (6,1) node[above,midway] {$z_r$};
\node (u1) at (8,3) {$x_1$};
\node at (8,2.15) {$\vdots$};
\node (um) at (8,1) {$x_m$};
\draw [<-, thick] (7,3) -- (u1);
\draw [<-, thick] (7,1) -- (um);
\end{tikzpicture}};
\end{center}
\caption{A multivariate nonlinear vector function $\vec{f}(\vec{x})$ can be represented in a decoupled representation $\vec{f}(\vec{x}) = \mat{W} \vec{g}(\mat{V}^\top \vec{x})$. A decoupled representation typically has fewer parameters, and reveals in an intuitive way, the internal nonlinearities of $\vec{f}(\vec{x})$.}
\label{fig:schematic}
\end{figure*}

This paper studies the question how a given nonlinear multivariate vector function $\vec{f} : \mathbb{R}^m \to \mathbb{R}^n$ can be decomposed into a simpler structure, as in~\citet{comon2017xrank,dreesen2014decoupling,tiels2013fctdpripwm,usevich2014mtns,vanmulders2014}. 
In particular, we investigate a structure of the form 
\begin{equation}
\label{eq:fWgVTx}
\vec{f}(\vec{x}) = \mat{W} \vec{g}(\mat{V}^\top \vec{x}),
\end{equation}
where $\mat{W}$ and $\mat{V}$ are transformation matrices, and the vector function $\vec{g}(\vec{z}) = \left[\begin{array}{ccc} g_1(z_1) & \cdots & g_r(z_r) \end{array} \right]^\top$ is composed of univariate functions $g_i(z_i)$ in its $r$ components.
The decoupled representation is visualized in Figure~\ref{fig:schematic}.

The decomposition~(\ref{eq:fWgVTx}) provides a representation that is easier to comprehend, as the nonlinearity is contained in a set of univariate components.
Moreover, it typically has a lower parametric complexity, and could hence be viewed as a form of nonlinear model order reduction. 

\subsection{Tensorization methods for decoupling polynomials}
When $\vec{f}$ is polynomial, the decomposition~(\ref{eq:fWgVTx}) has connections to the canonical polyadic decomposition (CPD) of a partially symmetric tensor (possibly a joint decomposition), see~\citet{usevich2017tensorizations}. 
Indeed, typical tensorization methods make use of the connection between homogeneous polynomials and their coefficient tensors~\citep{usevich2017tensorizations,comon2008symtensors}. 
The (homogeneous) scalar polynomial case, \ie $n=1$, is known as the Waring decomposition~\citep{carlini2003waringseveralforms,oeding2013eigtensors}, and is a fundamental problem in algebraic geometry. 
The problem~(\ref{eq:fWgVTx}) we study is hence very reminiscent of the classical Waring problem, but we consider the \emph{non-homogeneous case} of \emph{several polynomials}. 
The non-homogeneous Waring problem is studied in \cite{bialynicki2008,schinzel2002}.
The simultaneous Waring problem for several homogeneous polynomials is studied in \cite{carlini2003waringseveralforms,tiels2013fctdpripwm}. 

\subsection{Contributions and organization of this paper}
We start from the tensorization method of~\citet{dreesen2014decoupling}. 
 In this framework, the function $\vec{f}$ and its first-order derivative information are evaluated in a number of sampling points. 
A tensor is constructed from the set of corresponding Jacobian matrices, which admits a CPD that allows for the reconstruction of the decomposition~(\ref{eq:fWgVTx}).
This approach has the following advantages: (i) the order and size of the constructed tensor do not increase with the degree of $\vec{f}$, and, (ii) the approach is not limited to the use of polynomials.

In the current paper, we generalize the method~\citet{dreesen2014decoupling} to incorporate second-order derivative information. 
Since second-order derivative information leads to partially symmetric tensors, we are ultimately able to formulate the decomposition as a partially symmetric joint tensor decomposition. 
By involving the second-order derivatives, we impose additional constraints on the (joint) tensor decompositions, hence it is expected to enjoy more relaxed uniqueness conditions. 
In the article, we assume that an exact and uniquely identifiable~\cite{comon2017xrank} representation of $\vec{f}(\vec{x})$ exists. 
Nevertheless, the resulting joint tensor decomposition will ultimately be phrased as an optimization problem, and provides a natural starting point for studying both the noisy decoupling problem, as well as a model reduction interpretation, but this is beyond the scope of the current paper.  

The current article is organized as follows. 
Section~\ref{sec:method} outlines the tensor-based decoupling method, leading to a joint tensor decomposition formulation. 
We illustrate how the uniqueness properties improve by including second-order derivatives.
Section~\ref{sec:exp} validates the method on three simulation examples. 
Section~\ref{sec:concl} summarizes the results and points out a few future research directions.

\subsection{Notation}
Scalars are denoted by lowercase or uppercase letters and vectors are denoted by lowercase bold-face letters. 
Elements of a vector are denoted by lowercase letters with an index as subscript, \eg $\vec{x} = \left[ \begin{array}{ccc} x_1 & \ldots & x_m \end{array} \right]^\top$.
Matrices are denoted by uppercase bold-face letters, \eg $\mat{V} \in \mathbb{R}^{m \times r}$. 
The entry in the $i$-th row and $j$-th column of the matrix $\mat{V}$ is denoted by $v_{ij}$. 
A matrix $\mat{V} \in \mathbb{R}^{m \times r}$ has columns $\vec{v}_i$ as in $\mat{V} = \left[ \begin{array}{ccc} \vec{v}_1 & \ldots & \vec{v}_r \end{array} \right]$.
The transpose of a matrix $\mat{V}$ is denoted by $\mat{V}^\top$. 
A diagonal matrix with diagonal elements $a_1$, $a_2$, $a_3$ is denoted by $\operatorname{diag}(a_1, a_2, a_3)$ or $\operatorname{diag}(a_i)$.
Higher-order tensors are denoted by bold-face uppercase caligraphical letters, \eg $\ten{J} \in \mathbb{R}^{n \times m \times N}$. 
For scalar, vector, matrix and higher-order tensor \textit{functions}, we employ the same conventions.  
The outer product is denoted by $\circ$ and defined as follows: For $\ten{X} = \vec{u} \circ \vec{v} \circ \vec{w}$, the entry in position $(i,j,k)$ is $u_i v_j w_k$.
The Frobenius norm of a tensor $\ten{X}$ is denoted as $\left\| \ten{X} \right\|_F$.
The Euclidean norm of a vector $\vec{x}$ is denoted as $\left\| \vec{x} \right\|$.
The first-order and second-order derivatives of a univariate function $g(z)$ are denoted by $g'(z)$ and $g''(z)$, respectively.

\section{Decoupling multivariate functions using tensors}\label{sec:method}

\subsection{The Canonical Polyadic Decomposition}
The canonical polyadic decomposition (CPD)~\citep{carroll1970,harshman1970,kolda2009tdaa} is the decomposition of a tensor into a minimal sum of rank-one components.
For instance, a third-order tensor $\ten{T}$ has a CPD of the form
\begin{equation}
\ten{T} = \sum_{i=1}^R \vec{a}_i \circ \vec{b}_i \circ \vec{c},
\end{equation}
or in a short-hand notation $\ten{T} = \llbrack \mat{A}, \mat{B}, \mat{C} \rrbrack$, where $\mat{A} = \left[ \begin{array}{ccc} \vec{a}_1 & \cdots & \vec{a}_R \end{array} \right]$ (similar for $\mat{B}$ and $\mat{C}$).
The CPD is a celebrated tensor decomposition, which has found a variety of applications in signal processing and data sciences. 
One of the attractive properties of the CPD is its more relaxed uniqueness conditions. 
In contrast with matrix factorization, where uniqueness is only possible by imposing additional constraints (\eg orthogonality in the singular value decomposition (SVD)), the CPD has milder uniqueness properties~\cite{kruskal1977,delathauwer2006albtcdimaasmd,kolda2009tdaa,domanov2013uniqcpd1,domanov2013uniqcpd2}.

In our framework, the uniqueness conditions of the CPD will allow us to ensure that the proposed decoupling method retrieves the uniquely identifiable model structure: for identifiable models, uniqueness of the CPD is a sufficient condition for uniqueness of the model.  

\subsection{Decoupling functions using first-order information}
\begin{figure}[!htb]
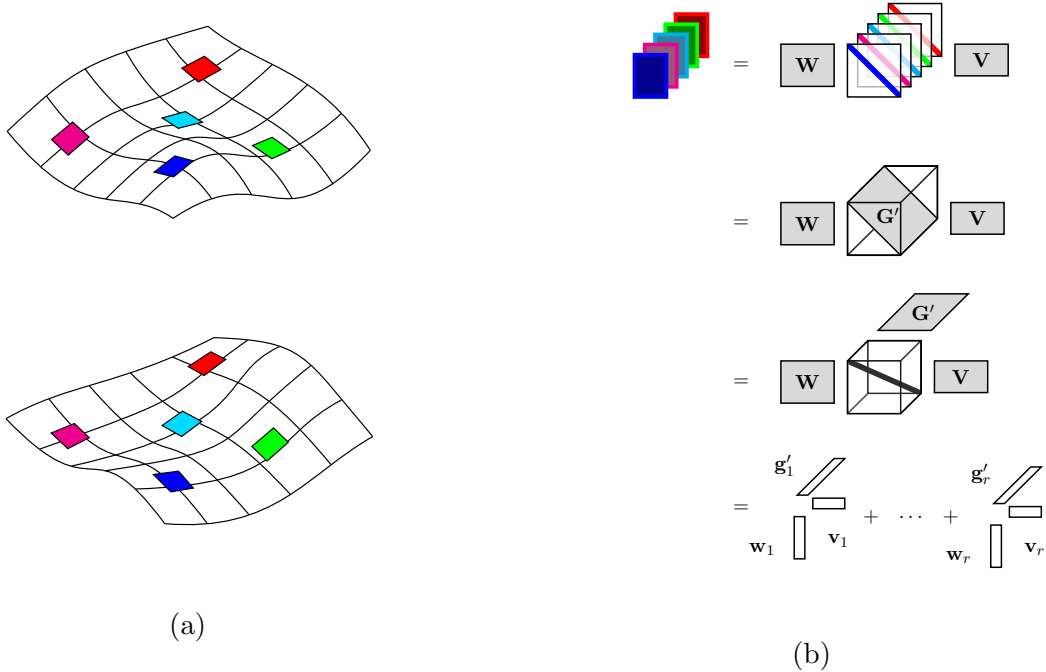

\begin{minipage}[c]{0.45\textwidth}
\centering
\begin{tikzpicture}[shift={(0,5cm)},anchor=center] 
    \node[scale=0.3]{\input{fig-twofunctions.tex} };
\end{tikzpicture}
\begin{center}
\vspace{0.5cm}
(a)
\end{center}
\end{minipage}
\begin{minipage}[c]{0.5\textwidth}
\centering
\begin{tikzpicture}[anchor=center,baseline] 
    \node[scale=0.70]{\input{fig-jaccpd.tex}};
\end{tikzpicture}
\begin{center}
\vspace{0.5cm}
(b)
\end{center}
\end{minipage}
\caption{The first-order information of $\vec{f}$ is collected in a set of sampling points $\vec{x}^{(k)}$, with $k=1,\ldots,N$ (indicated by the colored patches on the surfaces in (a)). 
The corresponding Jacobian matrices $\mat{J}(\vec{x}^{(k)})$ are arranged into a three-way tensor (b).
Each Jacobian matrix can be written as $\mat{J}(\vec{x}^{(k)}) = \mat{W} \operatorname{diag}(g_i'(\vec{v}_i^\top \vec{x}^{(k)})) \mat{V}^\top$. 
This results in a simultaneous matrix diagonalization problem, which is computed by the CPD.}
\label{fig:jaccpd}
\end{figure}

Consider the function $\vec{f}: \mathbb{R}^m \to \mathbb{R}^n$ that admits a decoupled representation~(\ref{eq:fWgVTx}). 
The Jacobian $\mat{J}(\vec{x})$ is represented with an $n \times m$ matrix function defined as
\begin{equation}
\mat{J}_{ij} (x) = \frac{\partial f_i (\vec{x})}{\partial x_j}. 
\end{equation}
The chain rule for derivation shows that the Jacobian matrix $\mat{J}(\vec{x})$ can be factorized as $\mat{J}(\vec{x}) = \mat{W} \operatorname{diag}( g_i'(\vec{v}_i^\top \vec{x}) ) \mat{V}^\top$,
or alternatively, in a CPD formulation 
$\ten{J} = \llbrack \mat{W}, \mat{V}, ( \vec{g}'(\mat{V}^\top \vec{x}) )^\top \rrbrack$.
Then an $n \times m \times N$ tensor $\ten{J}$, built from evaluating the Jacobian matrix $\mat{J}(\vec{x}^{(k)})$ in a set of $N$ sampling points $\vec{x}^{(k)}, k=1,\ldots,N$, admits a CPD of the form
\begin{equation}
\ten{J} = \llbrack \mat{W}, \mat{V}, \mat{G}' \rrbrack,
\label{eq:jaccpd}
\end{equation}
where $\mat{G}'$ contains the first-order derivatives of the functions $g_i$ in the $N$ points, \ie $\mat{G}'_{ki} = g_i'(\vec{v}_i^\top \vec{x}^{(k)})$.
In this way, the decoupled representation can be reconstructed from a simultaneous matrix diagonalization, or a CPD (Figure~\ref{fig:jaccpd}).

\subsection{Decoupling functions using second-order information}
Along the same lines, we may represent the Hessian $\ten{H}(\vec{x})$ by an $n \times m \times m$ tensor function defined as 
\begin{equation}
\ten{H}_{ijk}(\vec{x}) = \frac{\partial^2 f_i (\vec{x})}{\partial x_j \partial x_k}, 
\end{equation}
which is symmetric in the second and third mode since $\frac{\partial^2 f}{\partial x_i \partial x_j} =\frac{\partial^2 f}{ \partial x_j \partial x_i}$.
The Hessian tensor function has a CPD representation of the form
\begin{equation}
\ten{H}(\vec{x}) = \llbrack \mat{W}, \mat{V}, \mat{V}, (\vec{g}''(\mat{V}^\top \vec{x}))^\top\rrbrack. 
\end{equation}
By evaluating the Hessian in a set of $N$ sampling points $\vec{x}^{(k)}, k=1,\ldots,N$, we find the CPD of the $n \times m \times m \times N$ tensor $\ten{H}$ as 
\begin{equation}
\ten{H} = \llbrack \mat{W}, \mat{V}, \mat{V}, \mat{G}'' \rrbrack,
\label{eq:hesscpd}
\end{equation}
where the columns of $\mat{G}''$ contain the second-order derivatives of the functions $g_i$ in the $N$ points, \ie $\mat{G}''_{ki} = g_i''(\vec{v}_i^\top \vec{x}^{(k)})$.

%
%

\subsection{A joint tensor decomposition with partial symmetry}
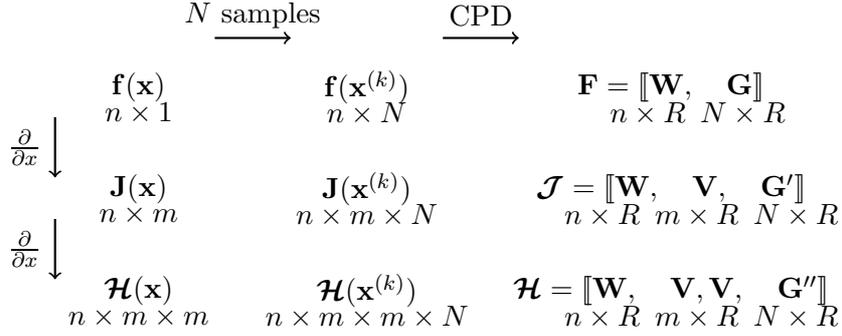
\begin{figure*}[!htb]
\begin{center}
\tikz \node[scale=1]{\begin{tikzpicture}

\node at (1.5, 0 + 0.35) {$\ten{H}(\vec{x})$};
\node at (1.5,0) {$n \times m \times m$};

\node at (4.5, 0 + 0.35) {$\ten{H}(\vec{x}^{(k)})$};
\node at (4.5, 0) {$n \times m \times m \times N$};

\node at (8.5, 0.35) {$\ten{H} = \llbrack \mat{W}, \quad \mat{V}, \mat{V},  \quad \mat{G''} \rrbrack$};
\node at (8.9, 0) {$n \times R \enspace m \times R \enspace N \times R$};


\node at (0, 0.9) {$\frac{\partial}{\partial x}$};
\draw [->, thick]  (0.4,1.3)  -- (0.4,0.5);


\node at (1.5, 0.35 + 1.0 + 0.35) {$\mat{J}(\vec{x})$};
\node at (1.5, 0.35 + 1.0) {$n \times m$};

\node at (4.5, 0.35 + 1.0 + 0.35) {$\mat{J}(\vec{x}^{(k)})$};
\node at (4.5, 0.35 + 1.0) {$n \times m \times N$};

\node at (8.5, 0.35 + 1.0 + 0.35) {$\ten{J} = \llbrack \mat{W}, \quad \mat{V},  \quad \mat{G'} \rrbrack$};
\node at (8.9, 0.35 + 1.0) {$n \times R \enspace m \times R \enspace N \times R$};


\node at (0, 2.25) {$\frac{\partial}{\partial x}$};
\draw [->, thick]  (0.4, 2.65)  -- (0.4, 1.85);


\node at (1.5, 0.35 + 1.0 + 0.35 + 1.0 + 0.35) {$\vec{f}(\vec{x})$};
\node at (1.5, 0.35 + 1.0 + 0.35 + 1.0) {$n \times 1$};

\node at (4.5, 0.35 + 1.0 + 0.35 + 1.0 + 0.35) {$\vec{f}(\vec{x}^{(k)})$};
\node at (4.5, 0.35 + 1.0 + 0.35 + 1.0) {$n \times N$};

\node at (8.5, 0.35 + 1.0 + 0.35 + 1.0 + 0.35) {$\mat{F} = \llbrack \mat{W},  \quad \mat{G} \rrbrack$};
\node at (8.85, 0.35 + 1.0 + 0.35 + 1.0) {$n \times R \enspace N \times R$};


\draw [->, thick] (2.5, 0.35 + 1.0 + 0.35 + 1.0 + 0.35 + 0.65) -- (3.5, 0.35 + 1.0 + 0.35 + 1.0 + 0.35 + 0.65);
\node at (3, 0.35 + 1.0 + 0.35 + 1.0 + 0.35 + 0.65 + 0.3) {$N$ samples};

\draw [->, thick] (5.5, 0.35 + 1.0 + 0.35 + 1.0 + 0.35 + 0.65) -- (6.5, 0.35 + 1.0 + 0.35 + 1.0 + 0.35 + 0.65);
\node at (6, 0.35 + 1.0 + 0.35 + 1.0 + 0.35 + 0.65 + 0.3) {CPD};

\end{tikzpicture}};
\end{center}
\caption{The proposed method takes into account a combination of first-order and second-order derivatives evaluations. These evaluations are organized into two tensors which admit a joint canonical polyadic decomposition with partial symmetry. The corresponding cost function is 
$\alpha_1\left\| \ten{J} - \llbrack \mat{W}, \mat{V}, \mat{G'} \rrbrack \right\|_F^2 + \alpha_2\left\| \ten{H} - \llbrack \mat{W}, \mat{V}, \mat{V}, \mat{G''} \rrbrack \right\|_F^2$.
}
\label{fig:jointdec}
\end{figure*}

The first-order and second-order derivative information can be combined into a joint tensor decomposition with partial symmetry. 
This can be phrased into the Structured Data Fusion framework~\citep{sorber2015sdf} and is implemented in tensorlab~\cite{tensorlab3} for MATLAB. 
The underlying optimization problem is 
\begin{equation}
\operatorname*{minimize}_{\mat{W}, \mat{V}, \mat{G'}, \mat{G''}} 
\alpha_1 \left\| \ten{J} - \llbrack \mat{W}, \mat{V}, \mat{G'} \rrbrack \right\|_F^2 + \alpha_2 \left\| \ten{H} - \llbrack \mat{W}, \mat{V}, \mat{V}, \mat{G''} \rrbrack \right\|_F^2,
\label{eq:optim}
\end{equation}
where the two terms in the cost function can be given different weights $\alpha_1$ and $\alpha_2$. 
The factor matrices $\mat{W}$ and $\mat{V}$ are shared among both decompositions. 
The partial symmetry in the Hessian tensor can be recognized in the fact that the factor $\mat{V}$ occurs twice. 
The joint decomposition approach is visualized in Figure~\ref{fig:jointdec}.

\subsection{Some remarks on uniqueness}
Our framework contains two notions that relate to `uniqueness' or `identifiability', which could be confusing. 
Therefore, it is useful to elaborate briefly on this matter. 
There is a notion of uniqueness at the level of the function decomposition~(\ref{eq:fWgVTx}), as well as at the level of the CPD of a (corresponding) tensor.
It is important to realize that in both cases, we consider so-called `essential uniqueness', which makes abstraction of trivial scaling and permutation invariances. 

On the one hand, in the polynomial case of~(\ref{eq:fWgVTx}), the problem at hand has a rich algebraic structure, which has recently been studied in the X-rank framework, leading to novel identifiability results~\cite{comon2017xrank}. 
These results assert that, for certain choices for $m$, $n$ and $r$, the function~(\ref{eq:fWgVTx}) has a single `unique' representation for generic choices of $\mat{W}$, $\mat{V}$ and $\vec{g}(\vec{z})$.
In other words, there does not exist an equivalent representation~(\ref{eq:fWgVTx}) having different $\mat{W}$, $\mat{V}$ and $\vec{g}(\vec{z})$. 

On the other hand, tensor decompositions have uniqueness properties themselves~\cite{kruskal1977,delathauwer2006albtcdimaasmd,kolda2009tdaa,domanov2013uniqcpd1,domanov2013uniqcpd2}, which ensure that if the tensor decomposition has converged to a (numerical) zero error, then the (essentially) unique underlying factorization has been retrieved.
Observe that uniqueness of the tensor decomposition is a sufficient condition for uniqueness of the function decomposition~(\ref{eq:fWgVTx}). 

Uniqueness properties of joint CPDs have only been studied recently in~\citep{sorensen2015ccpd1,sorensen2015ccpd2}.
Intuitively, it can be expected that, by imposing additional constraints on a decomposition, it is likely that uniqueness conditions are more easily met. 
 
\section{Numerical examples}\label{sec:exp}
In the current section, we will illustrate the proposed method on a number of examples. 
The joint decomposition~(\ref{eq:optim}) is implemented in tensorlab~\citep{tensorlab3}.

\subsection{Second-order derivatives for tackling the Waring decomposition}
In the single-output case, the Jacobian of the multivariate scalar function $f(\vec{x})$ is a vector function rather than a matrix function.
Its representation simplifies to $f(\vec{x}) = \vec{w}^\top \vec{g}(\mat{V}^\top \vec{x})$. 
Notice that $\vec{w}$ could be absorbed into the function $\vec{g}$, but we have chosen to keep it explicitly in the formula to show the resemblance to the vector function case.
For a scalar function $f(\vec{x})$, the Jacobian reduces to a vector function, rather than a matrix function, \ie $\vec{j}^\top(\vec{x}) = \llbrack \vec{w}^\top, \mat{V}, (\vec{g}'(\mat{V}^\top \vec{x})^\top \rrbrack$. 
Evaluating the Jacobian in a set of sampling points then gives rise to a matrix, rather than a tensor.
Summarizing, the $n \times m \times N$ Jacobian tensor $\ten{J}$ reduces in this situation to an $m \times N$ matrix $\mat{J}$, and we obtain a matrix factorization question, rather than the third-order tensor CPD~(\ref{eq:jaccpd}).
It is easy to understand that there is no unique solution, since one can insert $\mat{M} \mat{M}^{-1}$, with $\mat{M}$ an invertible $R \times R$ matrix, to obtain an equivalent factorization $\mat{J} = (\mat{V} \mat{M}) (\mat{G}' \mat{M}^{-T}) = \widetilde{\mat{V}} \widetilde{\mat{G}}'^\top$.

A possible solution to resolve this lack of uniqueness is to consider second-order derivative information of $f(\vec{x})$. 
We evaluate the $m \times m$ Hessian function 
\begin{equation}
\mat{H}_{ij}(\vec{x}) = \frac{\partial^2 f(\vec{x})}{\partial x_i \partial x_j},
\end{equation} 
in a set of $N$ sampling points $\vec{x}^{(k)}$. 
This gives an $m \times m \times N$ tensor $\ten{H}$ with
\begin{equation}
\ten{H} = \llbrack \vec{w}^\top, \mat{V}, \mat{V}, \mat{G}'' \rrbrack,
\end{equation}
as in~(\ref{eq:hesscpd}).
Now, we have again a case in which uniqueness of the CPD is attainable. 

For instance, we consider the function 
\begin{equation}
f(x_1,x_2) = - 37 x_1^3 - 213 x_1^2 x_2 - 399 x_1 x_2^2 + 5 x_1 - 239 x_2^3 + 9 x_2 - 2,
\end{equation}
which can be decomposed with 
\begin{equation}
\vec{w}^\top = \left[ \begin{array}{cc} 1& 1 \end{array} \right], \quad \mat{V} = \left[ \begin{array}{rr} 1 & 2 \\ 3 & 4 \end{array} \right],
\end{equation} 
and
\begin{equation}
\begin{array}{rcl} g_1(z_1) &=&  3 z_1^3 - z_1 + 5, \\ g_2 (z_2) &=& - 5 z_2^3 + 3 z_2 - 7.
\end{array}
\end{equation}
We draw an i.i.d.~set of sampling points $\vec{x}^{(k)}, k=1,\ldots,200$ from a uniform distribution between $-10$ and $10$ in both components.
The Jacobian matrix is hence a $2 \times 200$ matrix $\mat{J}$, which admits a non-unique rank-two factorization.
If we take into account the Hessian information, the problem again becomes a uniquely defined tensor question. 
The CPD of $\ten{H}$ is computed using tensorlab~\citet{tensorlab3}, and retrieves up to a scaling and permutation invariance, the true factors $\mat{V}$ and $\mat{G}''$. 

\subsection{Second-order derivatives improve uniqueness of the CPD}
For a general vector function $\vec{f}(\vec{x})$, including the Hessian information can improve the uniqueness properties beyond the Jacobian tensor method. 
For instance, in the case $m = n = 2$, the bound by~\citet{dreesen2014decoupling} ensures uniqueness up to $r \leq 2$. 
It can be verified that the Jacobian-based CPD method is not able to retrieve the underlying model. 
However,~\citet{comon2017xrank} asserts that $r=3$ is still identifiable for polynomial models (with degree $d\geq 3$). 
Considering the second-order information then leads to a CPD that is generically unique. 

We consider the function $\vec{f}(\vec{x})$ which is defined as 
\begin{equation}
\begin{array}{rcl}
f_1(x_1,x_2) &=& 24 x_1^3 + 36 x_1^2 x_2 - 4 x_1^2 + 18 x_1 x_2^2 - 4 x_1 x_2 + 84 x_2^3 - x_2^2 - 6 x_2 + 7,\\
f_2(x_1,x_2) &=& - 43 x_1^3 - 72 x_1^2 x_2 + 8 x_1^2 - 36 x_1 x_2^2 + 8 x_1 x_2 - 3 x_1 \\ & & \quad \quad \quad \quad \quad \quad \quad \quad \quad \quad \quad \quad \quad \quad  \quad+ 75 x_2^3 + 2 x_2^2 - 6 x_2 - 1,
\end{array}
\end{equation}
and admits a representation with 
\begin{equation}
\mat{W} = \left[ \begin{array}{rrr} 1 & 0 & 1\\-2 & -1 & 1 \end{array} \right] \mbox{ and } \mat{V} = \left[ \begin{array}{rrr} 2 & 1 & 0\\1 & 0 & 3 \end{array} \right],
\end{equation}
and 
\begin{equation}
\begin{array}{rcl}
g_1(z_1) &=& 3 z_1^3 - z_1^2 + 5, \\
g_2(z_2) &=& - 5 z_2^3 + 3 z_2 - 7,\\
g_3(z_3) &=& 3 z_3^3 - 2 z_3 + 2.
\end{array}
\end{equation}
A set of sampling points $\vec{x}^{(k)}, k=1,\ldots,200$ is sampled again uniformly on $[-10, 10]^2$.
Applying the Jacobian method results in an $2 \times 2 \times 200$ tensor which does not satisfy the uniqueness conditions.
Indeed, although the CPD has an error of the order $10^{-11}$, it does not return the correct factors. 
However, if we compute the CPD of the Hessian tensor, the underlying representation is found.

\subsection{Can we go beyond identifiable structures?}
In our experiments, we have observed a number of cases where adding Hessian info ensures interpretability, while the underlying model structure does not seem to be identifiable. 
For instance consider the $m=n=2$ and $r=4$ case of polynomials of degree $d=3$. 
We consider a function $\vec{f}(\vec{x})$ of the form~(\ref{eq:fWgVTx}) with
\begin{equation}
\mat{W} = \left[ \begin{array}{rrrr} 1  &   0  &   1   &  2\\
    -2 &   -1   &  1   &  3 \end{array} \right],
    \mbox{ and } 
\mat{V} = \left[ \begin{array}{rrrr}  2  &   1  &   0  &   1 \\
     1  &   0  &   3  &  -1  \end{array}
 \right],
\end{equation}
and 
\begin{equation}
\begin{array}{rcl}
g_1(z_1) &=& 3 z_1^3 - z_1^2 + 5,\\
g_2(z_2) &=& - 5 z_2^3 + 3 z_2 - 7, \\
g_3(z_3) &=& 3 z_3^3 - 2 z_3 + 2,\\
g_4(z_4) &=& z_4^3 - 2 z_4^2 + 1.
\end{array}
\end{equation}

We observe that by considering the second-order derivatives information, we are able to retrieve a decomposition having (numerical) zero error. 
However, the linear transformations $\mat{W}$ and $\mat{V}$ are \emph{not equal} (up to scaling and permutation) to the underlying factors. 
Nevertheless, when investigating the factor $\mat{G}''$, we see that the retrieved factor does have in its components a set of linear relations as expected from $\mat{G}''_{ki} = g_i''(\vec{v}_i^\top \vec{x}^{(k)})$. 
The fact that an interpretable model is obtained is a surprising result: it seems to suggest that, considering the second-order information enforces that only `interpretable' models are retrieved. 
However, we should mention that we have observed that this effect does not always hold for other cases.



\section{Conclusions and perspectives}\label{sec:concl}
In this article, we generalized a tensor-based method for finding a decoupled representation of a given nonlinear multivariate vector function. 
The method works by evaluating second-order derivatives in a set of sampling points.
First-order and second-order information can be combined in this way into a simultaneous higher-order tensor decomposition task with partial symmetry. 
We illustrated the promising abilities of this approach on a number of simulation examples.
The method was shown to outperform the existing approach: uniquely identifiable structures are recovered in a greater number of configurations, such as the single-output case, and cases in which the Jacobian tensor method was not able to ensure uniqueness.
We also observed that the method seems to extract meaningful representations even in some cases when the model structure is not uniquely identifiable, which is a property that merits further investigation.  

In future work, we want to investigate how function evaluations can be taken into account: this seems to make sense when the corresponding (matrix) factorization is a low-rank approximation, which occurs when $r<n$.
Also the use of higher-order derivatives is a possible extension, which might improve the uniqueness conditions even further. 
In this sense, the results in this article are a step towards a more general data-driven and noise-robust tensor-based framework for decoupling function representations.

\section*{Acknowledgments} 
This work was supported in part by the Flemish Government (Methusalem), and by the Fonds Wetenschappelijk Onderzoek -- Vlaanderen under EOS Project no 30468160 and research projects G.0280.15N and G.0901.17N.

\bibliographystyle{plain}
\bibliography{references}

\end{document}